\documentclass[12pt]{article}
\usepackage{amsfonts}

\newtheorem{lm}{Lemma}
\newtheorem{thm}{Theorem}

\newcommand{\ov}{\overline }

\begin{document}

\title{There are no structurally stable diffeomorphisms of odd-dimensional manifolds with codimension one non-orientable expanding attractors}
\author{V.~Medvedev\and E.~Zhuzhoma}
\date{Dedicated to Carlos Gutierrez on his 60th birthday}
\maketitle

\begin{abstract}
We prove that a structurally stable diffeomorphism of closed $(2m+1)$-manifold, $m\ge 1$, has no codimension one non-orientable expanding attractors.
\end{abstract}

\section{Introduction}

Structurally stable diffeomorphisms exist on any closed manifold (say a diffeomorphism $f$ structurally stable if all diffeomorphisms $C^1$-close to $f$ are conjugate to $f$). It is natural to study the question of
existence of such diffeomorphisms with some additional conditions. The condition we consider here is the presence of a codimension one non-orientable expanding attractor. Due to well known example of Plykin \cite{Plykin74}, the answer is YES for 2-manifolds. Medvedev and Zhuzhoma \cite{MedvZh2002} proved that for 3-manifolds the answer is NO. In the paper, we generalize the result of \cite{MedvZh2002} proving that there are no structurally stable diffeomorphisms with a codimension one non-orientable expanding attractor on closed odd-dimensional manifolds. The proof is shorter than \cite{MedvZh2002} and includes $d = 3$. As to orientable attractors, the answer is YES for any $d\ge 2$. Namely, starting with a codimension one Anosov diffeomorphism of the $d$-torus $T^d$, $d\ge 2$, a structurally stable diffeomorphism of $T^d$  with an orientable codimension one expanding attractor can be obtained by Smale's surgery \cite{Smale67}, so-called $DA$-diffeomorphism (see also \cite{KatokHassenblat95}, \cite{Plykin84}, \cite{Robinson-book99}).

Before the formulation of exact result, we give necessary definitions and notions.
Let $ f: M\to M $ be a diffeomorphism of a closed $d$-manifold $M$, $d = \dim M\ge 2$, endowed with some Riemann metric $\rho$ (all definitions in this section can be found in \cite{KatokHassenblat95} and \cite{Robinson-book99}, unless otherwise indicated). A point $x\in M$ is {\it non-wandering} if for any neighborhood $U$ of $x$, $f^n(U)\cap U \neq \emptyset$ for infinitely many integers $n$. Then the non-wandering set $NW(f)$, defined as the set of all non-wandering points, is an $f$-invariant and closed. A closed invariant set $\Lambda \subset M$ is {\it hyperbolic} if there is a continuous $f$-invariant splitting of the tangent bundle $T_{\Lambda}M$ into stable and unstable bundles $E^s_{\Lambda}\oplus E^u_{\Lambda}$ with
$$ \Vert df^n(v)\Vert \le C\lambda ^n\Vert v\Vert ,\quad \Vert df^{-n}(w)\Vert \le C\lambda ^n\Vert w\Vert ,
\quad \forall v\in E^s_{\Lambda}, \forall w\in E^u_{\Lambda}, \forall n\in \mathbb{N},$$
for some fixed $C > 0$ and $\lambda < 1$.  For each $x\in \Lambda$, the sets
$ W^s(x) = \{y\in M: \lim_{j\to \infty}\rho (f^j(x),f^j(y))\to 0 $, 
$ W^u(x) = \{y\in M: \lim_{j\to \infty} \rho (f^{-j}(x),f^{-j}(y)) \to 0 $
are smooth, injective immersions of $E_x^s$ and $E_x^u$ that are tangent to $E_x^s$, $E_x^u$ respectively.  $W^s(x)$, $W^u(x)$ are called {\it stable} and {\it unstable manifolds} at $x$.

For a diffeomorphism $ f: M\to M $, Smale \cite{Smale67} introduced the Axiom A: $NW(f)$ is hyperbolic and the periodic points are dense in $NW(f)$. A diffeomorphism satisfying the Axiom A is called $A$-diffeomorphism. According to Spectral Decomposition Theorem, $NW(f)$ of an $A$-diffeomorphism $f$ is decomposed into finitely many disjoint so-called basic sets $B_1$, $\ldots , B_k$ such that each $B_i$ is closed, $f$-invariant and contains a dense orbit.

A basic set $\Omega$ is called an {\it expanding attractor} if there is a closed neighborhood $U$ of $\Omega$ such that $f(U)\subset int~U$, $\cap _{j\ge 0}f^j(U) = \Omega$, and the topological dimension $\dim \Omega$ of $\Omega$ is equal to the dimension $\dim (E^u_{\Omega})$ of the unstable splitting $ E^u_{\Omega} $ (the name is suggested in \cite{Williams67}, \cite{Williams74}). $\Omega$ is codimension one if $\dim~\Omega = \dim~M - 1$. It is well known that a codimension one expanding attractor consists of the $(d-1)$-dimensional unstable manifolds $W^u(x)$, $x\in \Omega$, and is locally homeomorphic to the product of $(d-1)$-dimensional Euclidean space and a Cantor set. $W^s(x)$ is homeomorphic to $\mathbb{R}$ and can be endowed with some orientation. $W^u(x)$ is homeomorphic to $\mathbb{R}^{d-1}$ and can be endowed with some normal orientation (even if $M$ is non-orientable). Due to hyperbolic structure, any $W^s(x)$ intersects $W^u(x)$ transversally, $x\in \Omega$. Following \cite {Grines75},  say that $\Omega$ is {\it orientable} if for every $x\in \Omega$ the index of the intersection $W^s(x)\cap W^u(x)$ does not depend on a point of this intersection (it is either $+1$ or $-1$). The main result is the following theorem.
\begin{thm}\label{main-theorem}
Let $f: M\to M$ be a structurally stable diffeomorphism of a closed $(2m+1)$-manifold $M$, $m\ge 1$. Then the spectral decomposition of $f$ does not contain codimension one non-orientable expanding attractors.
\end{thm}
Our proof does not work for even-dimensional manifolds for which the existence of codimension one non-orientable expanding attractors stay open question (except $d = 2$).

{\it Acknowledgment}. The research was partially supported by RFFI grant 02-01-00098. We thank Roman Plykin and Santiago Lopez de Medrano for useful discussions.
This work was done while the second author was visiting Rennes 1 University (IRMAR) in March-Mai 2004. He thanks the support CNRS which made this visit possible. He would like to thank Anton Zorich and Vadim Kaimanovich for their hospitality.

\section{Proof of the main theorem}

Later on, $\Omega$ is a codimension one non-orientable expanding attractor of diffeomorphism $f: M\to M$. A point $p\in \Omega$ is called {\it boundary} if at least one component of $W^s(p) - p$ does not intersect $\Omega$. Boundary points exist and satisfy to the following conditions \cite{Grines75}, \cite{Plykin74}:
\begin{itemize}
\item There are  finitely many boundary points and each is periodic.
\item Given a boundary point $p\in \Omega$, there is a unique component of $W^s(p) - p$ denoted by $W^s_{\emptyset}(p)$ which does not intersect $\Omega$.
\item Given a point $x\in W^u(p) - p$, there is a unique arc $(x,y)^s\subset W^s(x)$ denoted by $(x,y)^s_{\emptyset}$ such that $(x,y)^s\cap \Omega = \emptyset$ and $y\in \Omega$.
\end{itemize}
An unstable manifold $W^u(p)$ containing a boundary point is called a {\it boundary unstable manifold}. Due to \cite{Grines75} and \cite{Plykin84}, the accessible boundary of $M - \Omega$ from $M - \Omega$ is a finite union of boundary unstable manifolds that splits into so-called bunches defined as follows. The family 
$W^u(p_1)$, $\ldots , W^u(p_k)$ is said to be a {\it $k$-bunch} if there are points $x_i\in W^u(p_i)$ and arcs $(x_i,y_i)^s_{\emptyset}$, $y_i\in W^u(p_{i+1})$, $1\le i\le k$, where $p_{k+1} = p_1$, $y_k\in W^u(p_1)$, and there are no $(k+1)$-bunches containing the given one.
\begin{lm}\label{non-orientability}
Let $f: M\to M$ be an $A$-diffeomorphism of a closed $(2m+1)$-manifold $M$, $m\ge 1$. If the spectral decomposition of $f$ contains a codimension one non-orientable expanding attractor, then $M$ is non-orientable.
\end{lm}
{\it Proof}. The non-orientability of $ \Omega $ implies that $ \Omega $ has at least one 1-bunch, say $W^u(p)$ \cite{Plykin84}. Therefore, given any point $x\in W^u(p) - p$, there is a unique point $y\in W^u(p) - p$ such that $(x,y)^s = (x,y)^s_{\emptyset}$, and vise versa. Let the map $ \phi : W^u(p) - p \to W^u(p) - p $
be given by $\phi (x) = y$ whenever $(x,y)^s = (x,y)^s_{\emptyset}$. Then $\phi$ is an involution, $\phi ^2 = id$.

Let $r$ be the period of $p$. Since the stable (as well as unstable) manifolds are $f$-invariant,
$ f^r\circ \phi |_{W^u(p)} = \phi \circ f^r|_{W^u(p)}$. Since the restriction $f^r|_{W^u(p)}$ is an expantion map with the unique hyperbolic fixed point $p$, $\phi$ can be extended homeomorphically to $W^u(p)$ putting $\phi (p) = p$.
By theorem 2.7 and lemma 2.1 \cite{Plykin84}, $\phi$ is conjugate to the antipodal involution, i.e. there exist a homeomorphism $ h: W^u(p)\to \mathbb{R}^{d-1} $ (in the intrinsic topology of $W^u(p)$) and the involution $\theta : \mathbb{R}^{d-1}\to \mathbb{R}^{d-1}$ of the type $\vec v \to -\vec v$ such that 
$\theta \circ h = h\circ \phi$. This implies that there is the $(d-1)$-dimensional ball $B^{d-1}\subset W^u(p)$ such that $p\in B^{d-1}$, the boundary $\partial B^{d-1}\stackrel{\rm def}{=}S^{d-2}$ is tamely embedded in $W^u(p)$, and $S^{d-2}$ is $\phi$-invariant. Moreover, there is the annulus $ S^{d-2}\times [0,1]\subset W^u(p) $ foliated by $S^{d-2}_t = S^{d-2}\times \{t\}$, $t\in [0,1]$, $ S^{d-2} = S^{d-2}_0 $, such that every $ S^{d-2}_t $ is $\phi$-invariant and bounds the $(d-1)$-dimensional ball $B^{d-1}_t\subset W^u(p)$ containing $p$. Since 
$\phi ^2 = id$, the set 
$$ B^{d-1}_t\bigcup _{x\in S^{d-2}_t}[x, \phi (x)]\stackrel{\rm def}{=}P_t $$
is homeomorphic to the projective space $ \mathbb{R}P^{d-1} $ for every $t\in [0,1]$. Since $d - 1 = 2m$ is even, $P_t$ is non-orientable. For any $ x\in S^{d-2}_{t_1}$ and $ y\in S^{d-2}_{t_2}$ with $t_1\neq t_2$, 
$ [x, \phi (x)]^s_{\infty} \cap [y, \phi (y)]^s_{\infty} = \emptyset$. Hence the set
$$ \bigcup _{x\in S^{d-2}\times [0,1]}[x, \phi (x)] \subset M$$
is homeomorphic to $ \mathbb{R}P^{d-1}\times [0,1] $. Since $ \mathbb{R}P^{d-1}\times [0,1] $ is a non-orientable $d$-manifold, $M$ is non-orientable. $\Box$
\medskip

{\it Proof of theorem \ref{main-theorem}}. Assume the converse. Then the spectral decomposition of $f$ contains a codimension one non-orientable expanding attractor, say $\Omega$. According to lemma \ref{non-orientability}, $M$ is non-orientable. Let $ \ov M $ be an orientable manifold such that $ \pi : \ov M\to M $ is a (nonbranched) double covering for $M$. Then there exists a diffeomorphism $\ov f: \ov M\to \ov M$ that cover $f$, i.e., 
$f\circ \pi = \pi \circ \ov f$. It is easy to see that $\ov f$ is an $A$-diffeomorphism with a codimension one expanding attractor $\ov \Omega \subset \pi ^{-1}(\Omega )$. It follows from lemma \ref{non-orientability} and orientability of 
$\ov M$ that $ \ov \Omega $ is orientable.

Because of $f$ is a structurally stable diffeomorphism, $f$ satisfies to the strong transversality condition \cite{Mane88} which is a local condition. Since $\pi$ is a local diffeomorphism, $\ov f$ satisfies to the strong transversality condition as well. Hence, $\ov f$ is structurally stable \cite{Robinson76}.

Take a periodic point $p\in \Omega$ on the boundary unstable manifold $W^u(p)$ that is a 1-bunch. Then the preimage $\pi ^{-1}(W^u(p))$ is a 2-bunch of $\ov \Omega$ consisting of unstable manifolds $ W^u(p_1) $, 
$ W^u(p_2) $, where $\{p_1, p_2\} = \pi ^{-1}(p)$ are boundary periodic points of $\ov f$. It was proved in 
\cite{GrinesZh2000a}, \cite{GrinesZh2000b} that $W^s_{\emptyset}(p_1)$ and $W^s_{\emptyset}(p_2)$ belong to the unstable manifolds $W^u(\alpha _1)$ and $W^u(\alpha _2)$ respectively of the repelling periodic points 
$ \alpha _1 $, $\alpha ^{\prime}$ (possibly, $ \alpha _1 = \alpha ^{\prime} $). Moreover, there are repelling periodic points $ \alpha _1 $, $ \ldots , \alpha _{k+1} = \alpha ^{\prime} $ and saddle periodic points
$ P_1 = p_1 $, $ P_2 $, $ \ldots , P_{k+1} $, $ P_{k+2} = p_2 $, $ k\ge 0 $, of index $ d - 1 $ such that the following conditions hold:
\begin{enumerate}
\item The set
$$ l = P_1\cup W^s_{\emptyset}(P_1)\cup \alpha _1\cup W^s(P_2)\cup \ldots \cup \alpha _{k+1}\cup W^s_{\emptyset}(P_{k+2})\cup P_{k+2}$$
is homeomorphic to an arc with no self-intersections whose endpoints are $P_1$ and $P_{k+2}$.
\item $ l - (P_1\cup P_{k+2}) \subset \ov M - \ov \Omega $.
\item The repelling periodic points $\alpha _i$ alternate with saddle periodic points $P_i$ on $l$.
\end{enumerate}
It follows from $f\circ \pi = \pi \circ \ov f$ that $\pi$ maps the stable and unstable manifolds of $\ov f$ into the stable and unstable manifolds respectively of $f$. Since $ \pi (P_1) = \pi (P_{k+2}) = p $,
$$ \pi (W^s_{\emptyset}(P_1)) = \pi (W^s_{\emptyset}(P_2)), \quad \pi (\alpha _1) = \pi (\alpha _{k+1}).$$
Hence (if $k\ge 1$),
$$ \pi (W^s(P_2)) = \pi (W^s(P_{k+1})),\quad \pi (P_2) = \pi (P_{k+1}),
\quad \pi (\alpha _2) = \pi (\alpha _{k}),\quad \ldots .$$
Due to item (3) above, the number of all periodic points on $l$ equals $2k+3$ that is odd. As a consequence, there is either a periodic point $\alpha _i$ with $ \pi (W^s(P_i)) = \pi (W^s(P_{i+1})) $ or a periodic point $P_i$ with
$ \pi (W^s_1(P_i)) = \pi (W^s_2(P_{i})) $, where $ \pi (W^s_1(P_i)) $,  $\pi (W^s_2(P_{i})) $ are different components of $ W^s(P_i) - P_i $. In both cases, there is a point ($ \alpha _i$ or $P_i$) at which $\pi$ is not a local homeomorphism. This contradiction concludes the proof. $\Box$

\bigskip

Dept. of Diff. Equat., Inst. of Appl. Math. and Cyber., Nizhny Novgorod State University, Nizhny Novgorod, Russia

{\it E-mail address}: medvedev@uic.nnov.ru
\medskip

Dept. of Appl. Math., Nizhny Novgorod State Technical University, Nizhny Novgorod, Russia

{\it E-mail address}: zhuzhoma@mail.ru

{\it Current e-mail address}: zhuzhoma@maths.univ-rennes1.fr

\end{document}